\documentclass[a4paper,12pt]{article}

\usepackage{amssymb}
\newtheorem{lemma}{Lemma}

\title{A variation of Gronwall's lemma}

\author{Quang-Cuong Pham\\
LPPA, Coll\`ege de France \\ 
Paris, France\\
\texttt{cuong.pham@normalesup.org}}

\begin{document}

\maketitle

\begin{abstract}
  We prove a variation of Gronwall's lemma.
\end{abstract}

The formulation and proof of the classical Gronwall's lemma can be
found in \cite{Gik}. We prove here a variation of this lemma, which we
were not able to find in the literature. The main difference from
usual versions of Gronwall's lemma is that $-\lambda$ is
\emph{negative}.

\begin{lemma}
  \label{lemma:gronwall}
  Let $g:[0,\infty[ \to\mathbb{R}$ be a continuous function, $C$ a
  real number and $\lambda$ a \emph{positive} real number. Assume that
  \begin{equation}
    \label{eq:bound1}
    \forall u,t \quad 0\leq u \leq t  \quad  g(t)-g(u) \leq \int_u^t -\lambda
    g(s)+C ds
  \end{equation}
  Then
  \begin{equation}
    \label{eq:bound2}
    \forall t \geq 0 \quad g(t)\leq \frac{C}{\lambda}+
    \left[g(0)-\frac{C}{\lambda}\right]^+e^{-\lambda t}
  \end{equation}
  where $[\cdot]^+=\max(0,\cdot)$.
\end{lemma}

\paragraph{Proof}

\textbf{Case 1 :} $C=0$, $g(0)>0$.

Define $h(t)$ by
\[
\forall t \geq 0 \quad  h(t)=g(0)e^{-\lambda t}
\]
Remark that $h$ is positive with $h(0)=g(0)$, and satisfies
(\ref{eq:bound1}) where the inequality has been replaced by an
equality
\[
\forall u,t \quad 0\leq u \leq t  \quad h(t)-h(u) = -\int_u^t \lambda
h(s) ds
\]

Consider now the set $ S=\{t\geq 0 \ |\ g(t) > h(t) \} $.  If
$S=\emptyset$ then the lemma holds true. Assume by contradiction that
$S\neq \emptyset$. In this case, consider an element $a\in S$. One has
by definition $g(a)>h(a)$. Since $g(0)=h(0)$, one also has
$a>0$. Consider now 
\[
m=\inf \{ a'<a\ |\ \forall t \in ]a',a[ \quad  g(t)>h(t)\}
\]
By continuity of $g$ and $h$ and by the fact that $g(0)=h(0)$, one has
$g(m)=h(m)$. One thus also has $m<a$ and
\begin{equation}
  \label{eq:epsilon}
  \forall t\in ]m,a[ \quad g(t)>h(t)  
\end{equation}

Consider now $\phi(t)=g(m)-\lambda\int_m^t g(s)ds$. Equation 
(\ref{eq:bound1}) implies that
\[
\forall t\geq m \quad g(t)\leq \phi(t)
\]

In order to compare $\phi(t)$ and $h(t)$ for
$t\in ]m,a[$, let us differentiate the ratio $\phi(t)/h(t)$.
\[
\left(\frac{\phi}{h}\right)'=\frac{\phi'h-h'\phi}{h^2}
=\frac{-\lambda gh+\lambda h \phi}{h^2}=\frac{\lambda
  h(\phi-g)}{h^2}\geq 0
\]
Thus $\phi(t)/h(t)$ is increasing for $t\in ]m,a[$. Since
$\phi(m)/h(m)=1$, one can conclude that
\[
\forall t\in ]m,a[ \quad \phi(t)\geq h(t)
\]
which implies, by definition of $\phi$ and $h$, that
\begin{equation}
  \label{eq:t}
  \forall t\in ]m,a[ \quad \int_m^t g(s)ds \leq \int_m^t h(s)ds  
\end{equation}

Choose now a $t_0\in ]m,a[$. Then one has by 
(\ref{eq:epsilon}) 
\[
\int_m^{t_0}g(s)ds > \int_m^{t_0}h(s)ds
\]
which clearly contradicts (\ref{eq:t}).

\textbf{Case 2 :} $C=0$, $g(0)\leq0$

Consider the set 
$
S=\{t\geq 0 \ | \  g(t) > 0 \}
$.
If $S=\emptyset$ then the lemma holds true. Assume by contradiction
that $S\neq \emptyset$. In this case, consider an element $a\in S$. 
One has by definition $g(a)>0$. Since $g(0)\leq 0$, one also has
$a>0$. Consider now 
\[
m=\inf \{ a'<a\ |\ \forall t \in ]a',a[ \quad  g(t)>0\}
\]
By continuity of $g$ and by the fact that $g(0)\leq 0$, one has
$g(m)=0$. One thus also has $m<a$ and
\begin{equation}
  \label{eq:epsilon2}
  \forall t\in ]m,a[ \quad g(t)>0  
\end{equation}
Choose now a $t_0\in ]m,a[$. Equation (\ref{eq:bound1}) implies that
\[
g(t_0)\leq -\lambda \int_m^{t_0} g(s)ds \leq 0
\]
which clearly contradicts (\ref{eq:epsilon2}).

\textbf{Case 3 :} $C\neq 0$

Define $\hat{g}=g-C/\lambda$. One has
\[
\forall u,t \quad 0\leq u \leq t  \quad \hat{g}(t)-\hat{g}(u)=
g(t)-g(u)  \leq \int_u^t -\lambda g(s)+C ds = -\int_u^t \lambda \hat{g}(s) ds
\]
Thus $\hat{g}$ satisfies the conditions of Case 1 or Case 2, and as a
consequence
\[
\forall t\geq 0 \quad \hat{g}(t) \leq [\hat{g}(0)]^+e^{-\lambda t}
\]
The conclusion of the lemma follows by replacing $\hat{g}$ by
$g-C/\lambda$ in the above equation. $\Box$

\subsection*{Acknowledgments}

The author would like to thank N. Tabareau and J.-J. Slotine for their
helpful comments and V. Valmorin for having pointed out an error in an
earlier version of the manuscript.

\end{document}